\documentclass[11pt]{amsart}

\newtheorem{theorem}{Theorem}[section]
\newtheorem{lemma}[theorem]{Lemma}

\theoremstyle{definition}

\newtheorem{coro}[theorem]{Corollary}

\theoremstyle{remark}

\numberwithin{equation}{section}

\newcommand{\ve}{}
\newcommand{\Polygn}{{\Pi}}
\def\cW{{\mathcal W}}
\def\cM{{\mathcal M}}
\def\cB{{\mathcal B}}

\def\cG{{\mathcal G}}

\def\cP{{\mathcal P}}

\begin{document}

\title[Best Simultaneous Diophantine Approximations]{Best Simultaneous Diophantine Approximations under a
Constraint on the Denominator}



\author{Iskander Aliev}
\address{School of Mathematics, University of Edinburgh, James Clerk Maxwell
Building, King's Buildings, Mayfield Road, Edinburgh EH9 3JZ, UK.
Tel.: +44/131/650-5056 Fax: +44/131/650-6553} \curraddr{}
\email{I.Aliev@ed.ac.uk}
\thanks{}

\author{Peter Gruber}
\address{Institut f\"ur Diskrete Mathematik und Geometrie E104, Technische
Universit\"at Wien,  Wiedner Hauptstrasse 8-10, A-1040 Wien,
Austria. Tel.: +43/1/58801-104 06 Fax: +43/1/58801-104 96}
\curraddr{} \email{peter.gruber@tuwien.ac.at}



\thanks{}




\begin{abstract}
We investigate the problem of best simultaneous Diophantine
approximation under a constraint on the denominator, as proposed by
Jurkat. New lower estimates for optimal approximation constants are
given in terms of critical determinants of suitable star bodies.
Tools are results on simultaneous Diophantine approximation of
rationals by rationals with smaller denominator. Finally, the
approximation results are applied to the decomposition of integer
vectors.
\end{abstract}

\maketitle


\section{Introduction}

The first investigations of simultaneous Diophantine approximation
with constraints on the denominator are due to Jurkat \cite{Jurkat}.
Kratz \cite{Kratz2,Kratz} considered the following particular
problem: let ${x}\in\mathbb R^{k}$, $k\ge 2$, and
$g(\cdot)=||\cdot||_2$. As in Kratz  \cite{Kratz2}, define for $Q>0$
the {\em successive minima} $\lambda_i=\lambda_i({ x}, Q)$,
$i=1,\ldots,k+1,$ {\it of} ${x}$ {\it under the constraint} $|q|\le
Q$ as follows: $\lambda_i$ is the minimum of all $\lambda\ge 0$, for
which there are $i$ linearly independent vectors ${
p_j}=(p_{j1},\ldots,p_{jk},p_{jk+1})\in\mathbb Z^{k+1}$,
$j=1,\ldots,i$, such that
$$g(p_{j k+1}{ x}-(p_{j1},\ldots,p_{jk}))\le
\lambda \ \text{and}\ |p_{j k+1}|\le Q \ \text{for} \
j=1,\ldots,i.$$
It is known (see e. g. \cite{Kratz2}) that the product of the first
$k$ successive minima satisfies
$$\lambda_1\cdots\lambda_k=O\left (\frac{1}{Q}\right )\,.$$
In this paper we are interested in an optimal constant
$c=c(k,||\cdot||_2)$ such that
$$\lambda_1\cdots\lambda_k<\frac{c}{Q}.$$
Kratz proved in \cite{Kratz} that
$$c(2,||\cdot||_2)=\frac{2}{\sqrt{3}}.$$

Assume now, that $g(\cdot)$ is the distance function of a bounded
star body $K$ in $\mathbb R^{k}$. In the following we consider the
above problem for $g(\cdot)$ and show in Theorem \ref{Optimal} that
$$c(k,g)\ge\frac{1}{\Delta(K)},$$
where $\Delta(K)$ is the critical determinant of $K$. Let $\gamma_k$
denote the Hermite constant. Since the critical determinant of the
$k$-dimensional unit ball equals $\gamma_k^{-k/2}$, we conclude that
$$c(k,||\cdot||_2)\ge\gamma_k^\frac{k}{2}.$$

For recent results on the upper bound for $c(k,g)$, see \cite{AH}.

To obtain these results, we study in detail the simultaneous
approximation of rational numbers by rational numbers with smaller
denominator. Let ${n}=(n_1,\ldots , n_{k}, n_{k+1})\in \mathbb
Z^{k+1}$, $k\ge 2$, be an integer vector. Assume that $0<n_1\le
\ldots\le n_{k+1}$ and that $\gcd(n_1, \ldots, n_{k+1})=1$. Consider
the problem of approximating the rational vector
$(n_1/n_{k+1},\ldots,n_k/n_{k+1})$ by rational vectors of the form
$(m_1/m_{k+1},\ldots,m_k/m_{k+1})$ with $m_i\in\mathbb Z$,
$i=1,\ldots,k+1,$ and $0\le m_{k+1}<n_{k+1}$. More precisely, we
investigate the behavior of the points

\begin{itemize}
\item[\hbox to 1.1cm{\rm (1)\hfill}]
$\displaystyle\left(m_1-m_{k+1}\frac{n_1}{n_{k+1}},\ldots,m_k-m_{k+1}\frac{n_k}{n_{k+1}}\right)$\label{rational_points}
\end{itemize}
as ${m}=(m_1,\ldots,m_k,m_{k+1})$ ranges over $\mathbb Z^{k+1}$.
Since these points form a $k$--dimensional lattice $\Lambda({n})$
(see Section \ref{Lambda_n} for details), we make use of tools from
the geometry of numbers.
%
%
%
%
%
%
%

Given an arbitrary lattice $\Lambda\subset \mathbb Q^{k}$, we can
construct a sequence of integer vectors ${ n}(t)$ such that the
sequence of corresponding lattices $\Lambda({n}(t))$ after an
appropriate normalization tends to $\Lambda$.

\begin{theorem}
For any rational lattice $\Lambda$ with basis
$\{{b}_1,\ldots,{b}_k\}$, ${b}_i\in\mathbb Q^{k}$, $i=1,\ldots,k$
and for all rationals $\alpha_1, \ldots,\alpha_k$ with
$0<\alpha_1\le\alpha_2\le\cdots\le\alpha_k\le 1$, there exists an
arithmetic sequence ${\mathcal P}$ and a sequence
${n}(t)=(n_1(t),\ldots,n_k(t),n_{k+1}(t))\in\mathbb Z^{k+1}$,
$t\in\mathcal P$, such that
$$\gcd(n_1(t),\ldots,n_k(t),n_{k+1}(t))=1$$
and $\Lambda({n}(t))$ has a basis ${a}_1(t),\ldots,{a}_k(t)$ with

\begin{itemize}
\item[\hbox to 1.1cm{\rm (2)\hfill}]
$\displaystyle
a_{ij}(t)=\frac{b_{ij}}{d\,t}+O\left(\frac{1}{t^{2}}\right)\
\text{for}\ i,j=1,\ldots,k \label{asympt_aij},$
\end{itemize}
where $d\in\mathbb N$ is such that $d\, b_{ij}, d\,
\alpha_j\,b_{ij}\in\mathbb Z$ for all $i,j=1,\ldots,k$. Moreover,

\begin{itemize}
\item[\hbox to 1.1cm{\rm (3)\hfill}]
$\displaystyle
n_{k+1}(t)=\frac{d^{k}t^{k}}{\det\Lambda}+O(t^{k-1})\label{asympt_h}$
\end{itemize}
and

\begin{itemize}
\item[\hbox to 1.1cm{\rm (4)\hfill}]
$\displaystyle \alpha_i(t):=\frac{n_i(t)}{n_{k+1}(t)}
=\alpha_i+O\left(\frac{1}{t}\right). \label{asympt_alpha}$
\end{itemize}

\label{main_lemma}
\end{theorem}

Let $\alpha(K)$ denote the anomaly of a set $K$, and if $f$ is the
distance function of $K$, then both $\lambda_i(f,\Lambda)$ and
$\lambda_i(K,\Lambda)$ denote the $i$\,th successive minimum of the
lattice $\Lambda$ with respect to the set $K$.

\begin{theorem}\label{main_main}
Let $K$ be a bounded star body in $\mathbb R^{k}$ and let

$$\mathbb U^{k+1}=\{{ x}\in\mathbb Z^{k+1}: 0<x_1\le\cdots\le
x_{k+1},\,\, \gcd(x_1,\ldots, x_{k+1})=1\}.$$ Then

$$C(K):=\sup_{n\in\mathbb U^{k+1}} \frac{\lambda_1(K,\Lambda({n}))\cdots\lambda_k(K,\Lambda({n}))}
{\det\Lambda({n})}=\frac{\alpha(K)}{\Delta(K)}.$$ Moreover, for all
$\alpha_1, \ldots,\alpha_k\in\mathbb Q$ with
$0<\alpha_1\le\alpha_2\le\cdots\le\alpha_k\le 1$, there exists an
infinite sequence of integer vectors
${n}(t)=(n_1(t),\ldots,n_k(t),n_{k+1}(t))\in\mathbb U^{k+1}$, $t\in
T=\{t_1, t_2,\ldots\}$, such that

\begin{itemize}
\item[\hbox to 1.1cm{\rm (5)\hfill}]
$\displaystyle\lim_{{t\to\infty}\atop{t\in T}}
\frac{\lambda_1(K,\Lambda({
n}(t)))\cdots\lambda_k(K,\Lambda({n}(t)))}{\det\Lambda({n}(t))}=
C(K); \label{main_main_n_t}$
\end{itemize}
\begin{itemize}
\item[\hbox to 1.1cm{\rm (6)\hfill}]
$\displaystyle\lim_{{t\to\infty}\atop{t\in T}}
\frac{n_i(t)}{n_{k+1}(t)} =\alpha_i,\ i=1,\ldots,k
\label{main_main_alpha}.$
\end{itemize}
and

\begin{itemize}
\item[\hbox to 1.1cm{\rm (7)\hfill}]
$\displaystyle\lim_{{t\to\infty}\atop{t\in T}} n_{k+1}(t)= \infty
.\label{infinity_h}$
\end{itemize}
\end{theorem}

The proof of Theorem~\ref{main_main} is based on the following lemma
which is of independent interest. Let $D$ denote the set of
functions $f: \mathbb R^{k}\rightarrow [0,+\infty)$ which are
positive homogeneous of degree $1$. We also denote by $o$ the zero
vector.

\begin{lemma}
Let $\{f_t\}$ be a sequence of functions in $D$ which converges
uniformly on $||{x}||\le 1$ to a function $f$ in $D$ and let
$\{L_t\}$ be a sequence of lattices in $\mathbb R^{k}$ which
converges to a lattice $L$. Then

\begin{itemize}
\item[\rm{(i)}]
$\limsup\limits_{t\rightarrow\infty}\lambda_i(f_t,L_t)\le\lambda_i(f,L)$,
for $i=1,\ldots,k$;
\item[\rm{(ii)}] If, in addition, $f({ x})=0$ only
for ${x}={ o}$, then
$\lim\limits_{t\rightarrow\infty}\lambda_i(f_t,L_t)$ exists and
equals $\lambda_i(f,L)$ for $i=1,\ldots,k$.
\end{itemize}
\label{Super_Cool}
\end{lemma}
This lemma clearly implies the following result.

\begin{coro}
If $\{L_t\}$ is a sequence of lattices in $\mathbb R^{k}$ convergent
to a full lattice $L$ and $K$ is a bounded star body then
$$\lim_{t\rightarrow\infty}\lambda_i(K,L_t)=\lambda_i(K,L) \text{
for each } i=1,\dots,k.$$\label{Cool}
\end{coro}

A similar result about centrally symmetric convex bodies was
recently proved by the first author jointly with Schinzel and
Schmidt in \cite{We}.

\begin{theorem}
Let $g(\cdot)$ be the distance function of a bounded star body $K$
in $\mathbb R^{k}$. Then
$$c(k,g)\ge\frac{1}{\Delta(K)}.$$\label{Optimal}
\end{theorem}

Theorem~\ref{main_lemma} will be applied in
Section~\ref{section_of_decompositions} to the problem of
decomposition of integer vectors, where the problem is considered
with respect to the supremum norm. For recent results on this
problem for the Euclidean norm $||\cdot||_2$, see \cite{We}. By
tradition, we denote the supremum norm of a vector ${a}$ by
$h({a})$.

Given $m$ linearly independent vectors ${n}_1, \ldots, {n}_m$ in
$\mathbb Z^{k+1}$ let $H({ n}_1, \ldots, {n}_m)$ denote the maximum
of the absolute values of the $m\times m$--minors of the matrix
$({n}_1^{t}, \ldots, {n}_m^{t})$ and let $D({n}_1, \ldots, {n}_m)$
be the greatest common divisor of these minors. Then $h({n})=H({n})$
for ${n}\neq{o}$. For $k+1>l>m>0$ let

\begin{itemize}
\item[\hbox to 1.1cm{\rm (8)\hfill}]
$\displaystyle c_0(k+1,l,m)= \sup\inf \left(\frac{D({\ve n}_1,
\ldots, {\ve n}_m)}{H({\ve n}_1, \ldots, {\ve
n}_m)}\right)^{\frac{k-l+1}{k-m+1}}\prod_{i=1}^{l}h({\ve p}_i),
\label{c0}$
\end{itemize}
where the supremum is taken over all sets of linearly independent
vectors ${\ve n}_1, \ldots, {\ve n}_m$ in $\mathbb Z^{k+1}$ and the
infimum over all sets of linearly independent vectors ${\ve
p}_1,\ldots,{\ve p}_l$ in $\mathbb Z^{k+1}$ such that
$${\ve n}_i=\sum_{j=1}^{l}u_{ij}{\ve p}_j, u_{ij}\in \mathbb Q
\text{ for all } i\le m.$$ It has been proved in \cite{Sch4} that
for fixed $l,m$,

\begin{itemize}
\item[\hbox to 1.1cm{\rm (9)\hfill}]
$\displaystyle\limsup_{k\rightarrow\infty}c_0(k+1,l,m)<\infty
\label{lessinfty}$
\end{itemize}
and in \cite{DecompII} it was shown that
$$c_0(k+1,2,1)\le\frac{2}{(k+1)^\frac{1}{k}}.$$
A result in \cite{ChSch} says that $c_0(3,2,1)=2/\sqrt{3}$. Note
that $$c_0(k+1,2,1)=\sup_{{\ve n}\in\mathbb Z^{k+1}\setminus\{o\}}
\inf_{\begin{array}{c} \scriptstyle{\ve p},\scriptstyle{\ve
q}\in\mathbb Z^{k+1}\setminus\{o\}\\
\scriptstyle\dim({\ve p},{\ve q})=2\\
\scriptstyle{\ve n} = u{\ve p}+v{\ve q},u, v \in \mathbb Z
\end{array}
} \frac{h({\ve p})h({\ve q})}{h({\ve n})^{1-\frac{1}{k}}}.$$

In this paper we continue to study the behavior of $c_0(k+1,2,1)$
and prove the following theorem.

\begin{theorem}
For $k\ge 3$ $$\limsup_{
\begin{array}{c}
\scriptstyle{\ve n}\in\mathbb Z^{k+1}\\
\scriptstyle h({\ve n})\rightarrow\infty\\
\end{array}
} \inf_{
\begin{array}{c}
\scriptstyle{\ve p},\scriptstyle{\ve q}\in\mathbb Z^{k+1}\\
\scriptstyle\dim({\ve p},{\ve q})=2\\
\scriptstyle{\ve n} = u{\ve p}+v{\ve q},u, v \in \mathbb Z
\end{array}
} \frac{h({\ve p})h({\ve q})}{h({\ve n})^{1-\frac{1}{k}}}\ge
\frac{1}{(k+1)^{\frac{1}{k}}}.$$ \label{th_on_decomp}
\end{theorem}

A more general result of Cha{\l}adus \cite{one_half} yields a weaker
inequality with $1/2$ instead of $1/(k+1)^{1/k}$.

\section{The Lattice $\Lambda({\ve n})$, Rational Weyl
Sequences and Systems of Linear Congruences} \label{Lambda_n}

In this section we construct a special lattice $\Lambda({\ve n})$.
Its points correspond to points of the form (1). Given the vector
${\ve n}$, there is a basis of the lattice $\mathbb Z^{k+1}$ of the
form ${\ve n},{\ve v}_1,\ldots,{\ve v}_k$. Let
$${\ve v}'_i=\left(v_{i1}-v_{ik+1}\frac{n_1}{n_{k+1}},\ldots,v_{ik}-v_{ik+1}\frac{n_k}{n_{k+1}}\right)\in
\mathbb R^{k}, i=1,\ldots ,k.$$ The equality
$$A_1{\ve v}'_1+\ldots+A_k{\ve v}'_k={\ve o}$$
implies that
$$n_{k+1}A_1{\ve v}_1+\ldots+n_{k+1}A_k{\ve v}_k+A_{k+1}{\ve
n}={\ve o}$$ with $A_{k+1}=-A_1v_{1k+1}-\ldots-A_kv_{kk+1}$. Thus
the vectors ${\ve v}'_1,\ldots,{\ve v}'_k$ are linearly independent.
Denote by $\Lambda({\ve n})$ the $k$--dimensional lattice with basis
$\{{\ve v}'_1,\ldots,{\ve v}'_k\}$. Since
$$\begin{array}{l} 1=\det \left (\begin{array}{cccc}
v_{11} & \ldots & v_{k1} & n_1\\
\vdots & \ddots & \vdots & \vdots\\
v_{1k} & \ldots & v_{kk} & n_k\\
v_{1k+1} & \ldots & v_{kk+1} & n_{k+1}\\
\end{array}
\right)\\[8ex]
=n_{k+1}\det \left (
\begin{array}{cccc}
v_{11}-v_{1k+1}\frac{n_1}{n_{k+1}} & \ldots &
v_{k1}-v_{kk+1}\frac{n_1}{n_{k+1}} &
\frac{n_1}{n_{k+1}}\\
\vdots & \ddots & \vdots & \vdots\\
v_{1k}-v_{1k+1}\frac{n_k}{n_{k+1}} & \ldots &
v_{kk}-v_{kk+1}\frac{n_k}{n_{k+1}} &
\frac{n_k}{n_{k+1}}\\
0 & \ldots & 0 & 1\\
\end{array} \right ),
\end{array}
$$
we have $\det\Lambda({\ve n})=1/n_{k+1}$. It is easily seen, that
for every non--zero vector ${\ve v}\in\Lambda({\ve n})$ there is a
unique vector ${\ve m}\in\mathbb Z^{k+1}$ such that
$${\ve v}=\left(m_{1}-m_{k+1}\frac{n_1}{n_{k+1}},\ldots,m_{k}-m_{k+1}\frac{n_k}{n_{k+1}}\right),
\text{ where } 0\le m_{k+1}<n_{k+1}.$$ Thus there is a one--to--one
correspondence between the points of $\Lambda({\ve n}) \setminus
\{{\ve o}\}$ and the non--zero integer vectors with $0\le
m_{k+1}<n_{k+1}$. Note also that since ${\ve v}\neq {\ve o}$, the
vectors ${\ve m}$ and ${\ve n}$ are linearly independent.

If $\Lambda$ is a lattice, let $\Lambda^{*}$ be its {\it polar
lattice}, see \cite{GrLek}.
%
The lattice $\Lambda({ n})$ is related to the lattice
$\Lambda^{\bot}({ n})$ of integer vectors orthogonal to ${\ve n}$.
Let $\Lambda^{\bot}_{k+1}({ n})$ be the $k$-dimensional lattice
obtained by omitting the $(k+1)$st coordinate in $\Lambda^{\bot}({
n})$. Then the following holds:

\begin{lemma}
The lattice $\Lambda^{\bot}_{k+1}({ n})$ is the polar lattice of the
lattice $\Lambda({n})$,
$$\Lambda^{\bot}_{k+1}({n})=\Lambda({n})^{*}.$$\label{Pol}
\end{lemma}

The lattice $\Lambda({\ve n})$ appears in some problems of number
theory. Let $\theta_1,\ldots,\theta_k$, $k\ge 2$, be real numbers
and let $\cW_k$ be the sequence of $k$--dimensional vectors

\begin{itemize}
\item[\hbox to 1.1cm{\rm (10)\hfill}]
$(i\theta_1 \mod 1, \ldots , i\theta_k \mod 1), \ i=0,1,2\ldots
\label{Weyl_points}$
\end{itemize}

$\cW_k$ is called a $k$--dimensional {\it Weyl sequence}. We shall
consider the case where
$$\theta_1=\frac{n_1}{n_{k+1}},\ldots,\theta_k=\frac{n_k}{n_{k+1}}.$$
Then $\cW_k$ is $n_{k+1}$--periodic and the set
$$\Lambda(\cW_k)=\{{\ve x}+{\ve y}: {\ve x}\in\mathbb Z^{k},\ {\ve y}\in
\cW_k  \}$$ is a $k$--dimensional lattice. It can be shown easily
that
$$\Lambda(\cW_k)=\Lambda({\ve n}).$$

Consider the lattice $n_{k+1}\Lambda({\ve
n})=n_{k+1}\Lambda(\cW_k)\subset\mathbb Z^{k}$. The points in (10),
multiplied by $n_{k+1}$, can be written in the form
$$(i n_1 \mod n_{k+1}, \ldots , i n_k \mod n_{k+1}), \
i=0,1,2,\ldots$$ Therefore, any point $(x_1,\ldots,x_k)\in
n_{k+1}\Lambda({\ve n})$ is a solution of the system

\begin{itemize}
\item[\hbox to 1.1cm{\rm (11)\hfill}]
$\displaystyle \left \{
\begin{array}{cccc}
x_1+r n_1&   \equiv&   0&  (\mod n_{k+1})\,\\
\vdots\\
x_k+r n_k&   \equiv&   0&  (\mod n_{k+1})\,\\
\end{array}
\right. \label{congruences}$
\end{itemize}
where $r$ is an integer corresponding to $m_{k+1}$. Hence we may
consider Theorems \ref{main_lemma} and \ref{main_main} as results on
rational Weyl sequences and solutions of the system (11).

\section{Proof of Lemma \ref{Pol}}
\label{Polarity}

Let ${\ve v}$ be a primitive non--zero vector of $\Lambda({\ve n})$
and ${\ve V}=n_{k+1}{\ve v}$. Choose a vector ${\ve m}\in\mathbb
Z^{k+1}$ such that $${\ve
v}=\left(m_{1}-m_{k+1}\frac{n_1}{n_{k+1}},\ldots,m_{k}-m_{k+1}\frac{n_k}{n_{k+1}}\right)
.$$ Let $\Lambda({\ve m}, {\ve n})$ denote the lattice with basis
${\ve m}, {\ve n}$. Since ${\ve v}$ is primitive, we have that
$$\Lambda({\ve m}, {\ve n})= S({\ve m}, {\ve n})\cap\mathbb Z^{k+1},$$
where $S({\ve m}, {\ve n})$ denotes the subspace of $\mathbb
Q^{k+1}$ spanned by the vectors ${\ve m}, {\ve n}$.

Consider the lattice $\Lambda^{\bot}({\ve m}, {\ve n})$ of integer
vectors orthogonal to $S({\ve m}, {\ve n})$ and choose a basis

\begin{itemize}
\item[\hbox to 1.1cm{\rm (12)\hfill}]
$\begin{array}{llllllr} {\ve a}_1' & =  ( a_{11} , \ldots , a_{1k},
a_{1k+1})\,,\\
\vdots\\
{\ve a}_{k-1}' & =  ( a_{k-11} , \ldots ,  a_{k-1k},
a_{k-1k+1})\,,\\
{\ve a}_{k}' & =  ( a_{k1} , \ldots ,  a_{kk},
a_{kk+1})\,\\
\end{array}$
\label{star_basis}
\end{itemize}
of the lattice $\Lambda^{\bot}({\ve n})$ such that the first $k-1$
vectors ${\ve a}_1'\ldots {\ve a}_{k-1}'$ form a basis of
$\Lambda^{\bot}({\ve m}, {\ve n})$. It is easy to see that the
vectors
$$\begin{array}{lllll}
{\ve a}_1 & = & ( a_{11} &, \ldots , & a_{1k})\,,\\
\vdots\\
{\ve a}_{k} & = & ( a_{k1} &, \ldots , & a_{kk})\,\\
\end{array}$$
form a basis of $\Lambda^{\bot}_{k+1}({\ve n})$. Consider the matrix
$$A= \left (
\begin{array}{cccc}
a_{11} & \ldots  & a_{1k} & a_{1k+1}\\
&\ddots\\
a_{k-11} & \ldots  & a_{k-1k} & a_{k-1k+1}\\
\end{array}
\right )$$ and denote by $A_{ij}$ the minor obtained by omitting the
$i$th and $j$th columns in $A$.

Let
$${\ve V}'_i=m_i {\ve n} - n_i{\ve m}$$
and let ${\ve V}_i$ be the vector obtained by omitting the $i$th
entry in ${\ve V}'_i$ (note this entry is 0). When omitting the
$i$th entry, we preserve the numbering of the remaining entries. For
example, we consider ${\ve V}_3$ as a vector of the $k$--dimensional
space with coordinates $x_1$, $x_2$, $x_4$, \ldots, $x_{k+1}$. In
particular, ${\ve V}_{k+1}={\ve V}$. Let $\Lambda^{\bot}_i({\ve m},
{\ve n})$ denote the lattice obtained by omitting the $i$th entries
of all vectors of the lattice $\Lambda^{\bot}({\ve m}, {\ve n})$,
preserving the numbering of the remaining entries. Denote by
$V_{ij}$ the $j$th entry of ${\ve V}_i$. Then the following result
holds.

\begin{lemma}
$V_{ij}=\epsilon_{ij}A_{ij}$, where $\epsilon_{ij}=\pm 1$ and
$\epsilon_{k+1i}\epsilon_{k+1j}=(-1)^{i-j}$. \label{almost_polar}
\end{lemma}
\begin{proof}
${\ve V}'_i\in\Lambda({\ve m}, {\ve n})$ implies that ${\ve
V}'_i\bot\Lambda^{\bot}({\ve m}, {\ve n})$ and thus ${\ve
V}_i\bot\Lambda^{\bot}_i({\ve m}, {\ve n})$. Hence ${\ve V}_i$ can
be represented in the form

\begin{itemize}
\item[\hbox to 1.1cm{\rm (13)\hfill}]
${\ve V}_i=s_i (\text{external product of the vectors of a basis of
} \Lambda^{\bot}_i({\ve m}, {\ve n})),\ s_i\in\mathbb R.$
\label{extern_product}\end{itemize} Therefore,
$$V_{ij}=\epsilon_{ij} t_i A_{ij},\ \epsilon_{ij}=\pm 1,\ t_i>0$$
and clearly $\epsilon_{k+1i}\epsilon_{k+1j}=(-1)^{i-j}$. In order to
see this, it is enough to note that the basis $\{{\ve a}_1'\ldots
{\ve a}_{k-1}'\}$  of $\Lambda^{\bot}_i({\ve m}, {\ve n})$ obtained
from (12) is a basis of the lattice on the right hand side of (13).
Further, the equation $V_{ij}=-V_{ji}$ implies that $t_i=t_j$. Let
$t=t_1=\ldots=t_k$. It is well known that

\begin{itemize}
\item[\hbox to 1.1cm{\rm (14)\hfill}]
$\det\Lambda({\ve m}, {\ve n})=\det\Lambda^{\bot}({\ve m}, {\ve
n}),\label{dets}$
\end{itemize}
see e.g.~\cite{BombVaal}, p.~27/28. For the first determinant holds
$$\det\Lambda({\ve m}, {\ve n})=\left ( \begin{array}{cc}
{\ve m}{\ve m} & {\ve m}{\ve n}\\
{\ve m}{\ve n} & {\ve n}{\ve n}\\
\end{array}\right )=\frac{1}{2}\sum_{i\neq j}V_{ij}^{2}
=\frac{t^2}{2}\sum_{i\neq j}A_{ij}^{2}.$$ On the other hand, by the
Laplace identity (see e.g.~\cite{Schmidt}, Lemma 6D), we can write
the second determinant as
$$\det\Lambda^{\bot}({\ve m}, {\ve n})=\det({\ve
a}'_i{\ve a}'_j)_{i,j=1}^{k-1}=\frac{1}{2}\sum_{i\neq
j}A_{ij}^{2},$$ and by (14)  $t=t_1=\ldots=t_k=1$.
\end{proof}

Since ${\ve V}=n_{k+1}{\ve v}$, Lemma \ref{almost_polar} implies
that the vector ${\ve v}$ is orthogonal to the vectors ${\ve
a}_1,\ldots,{\ve a}_{k-1}$ and
$$\begin{array}{l}
{\ve v}{\ve a}_{k}=\frac{1}{n_{k+1}}{\ve V}{\ve
a}_{k}=\frac{1}{n_{k+1}}(V_{k+11}a_{k1}+\ldots+V_{k+1k}a_{kk})\\
=\pm\frac{1}{n_{k+1}}(A_{k+11}a_{k1}-A_{k+12}a_{k2}+\ldots+(-1)^{k-1}A_{k+1k}a_{kk})\\
=\pm\frac{1}{n_{k+1}}\det \Lambda_{k+1}^{\bot}({\ve n})=\pm 1
.\end{array}$$ By taking, if necessary, $-{\ve v}$ instead of ${\ve
v}$, we may assume that ${\ve v}\,{\ve a}_{k}=1$. This shows that
${\ve v}\in \Lambda_{k+1}^{\bot}({\ve n})^{*}$. Thus $\Lambda({\ve
n})$ is a sublattice of $\Lambda_{k+1}^{\bot}({\ve n})^{*}$. Since
$$\det\Lambda({\ve n})=\det(\Lambda_{k+1}^{\bot}({\ve
n}))^{*}=\frac{1}{n_{k+1}},$$ these lattices coincide.

\section{Proof of Theorem \ref{main_lemma}}

Let $\{{\ve b}^{*}_1,\ldots,{\ve b}^{*}_k\}$ be the basis of the
polar lattice $\Lambda^{*}$ given by 
$${\ve b}^{*}_i {\ve
b}_j= \left \{
\begin{array}{ll}
1, & i=j\,,\\
0, & \mbox{otherwise}\,.\\
\end{array}
\right.$$ We shall apply Theorem 1 of \cite{NewSL}, where $m=1$,
$F=1$, and $F_{1\nu}$, $\nu=1,\ldots,k+1$ are the minors of order
$k$ of the matrix
$$\begin{array}{l}
M=M(T,T_1,\ldots,T_k)\\
=\left (
\begin{array}{ccccc}
db^{*}_{11}T+T_1 & db^{*}_{12}T  & \ldots & db^{*}_{1k}T &
d\sum_{i=1}^{k}\alpha_ib^{*}_{1i}T\\
db^{*}_{21}T & db^{*}_{22}T+T_2  & \ldots & db^{*}_{2k}T &
d\sum_{i=1}^{k}\alpha_ib^{*}_{2i}T\\
\vdots & \vdots & & \vdots & \vdots\\
db^{*}_{k1}T & db^{*}_{k2}T  & \ldots & db^{*}_{kk}T+T_k &
d\sum_{i=1}^{k}\alpha_ib^{*}_{ki}T\\
\end{array}
\right ) , \end{array}$$ where $T, T_1, \ldots, T_k$ are variables.
Let $M_i=M_i(T,T_1,\ldots,T_k)$ and let $B^{*}_i$ be the minor
obtained by omitting the $i$th column in $M$ or in the matrix
$$\left (
\begin{array}{ccccc}
b^{*}_{11} & b^{*}_{12}  & \ldots & b^{*}_{1k} &
\sum_{i=1}^{k}\alpha_ib^{*}_{1i}\\
b^{*}_{21} & b^{*}_{22}  & \ldots & b^{*}_{2k} &
\sum_{i=1}^{k}\alpha_ib^{*}_{2i}\\
\vdots & \vdots & & \vdots & \vdots\\
b^{*}_{k1} & b^{*}_{k2}  & \ldots & b^{*}_{kk} &
\sum_{i=1}^{k}\alpha_ib^{*}_{ki}\\
\end{array}
\right ) ,$$ respectively. As in the proof of Theorem 2 in
\cite{NewSL} we have that

\begin{itemize}
\item[\hbox to 1.1cm{\rm (15)\hfill}]
$|B^{*}_{k+1}|=|\det(b^{*}_{ij})|\neq 0,$\label{last_coordinate}
\end{itemize}

\begin{itemize}
\item[\hbox to 1.1cm{\rm (16)\hfill}]
$|B^{*}_i|=\alpha_i|B^{*}_{k+1}|,$\label{coordinates}
\end{itemize}

\begin{itemize}
\item[\hbox to 1.1cm{\rm (17)\hfill}]
$M_i=d^{k}B^{*}_iT^{k}+\mbox{polynomial of degree less than}\, k\,
\mbox{in} T $\label{asympt_minors}
\end{itemize}
and $M_1,\ldots,M_k$ have no common factor. By Theorem~1 of
\cite{NewSL} there exist integers $t_1,\ldots,t_k$ and an arithmetic
progression $\cP$ such that, for $t\in\cP$, we have
$$\gcd(M_1(t,t_1,\ldots,t_k), \ldots ,M_{k+1}(t,t_1,\ldots,t_k))=1
.$$ Let $${\ve n}(t)=(M_1(t,t_1,\ldots,t_k), \ldots
,(-1)^{k}M_{k+1}(t,t_1,\ldots,t_k)).$$ Then (3)and (4) hold.

To prove the equality (2), consider the lattice
$\Lambda_{k+1}^{\bot}({\ve n}(t))$, $t\in\cP$ with basis
$$\begin{array}{llllllllllll}
{\ve a}_1^{*}(t) & = ( db^{*}_{11}t+t_1,  db^{*}_{12}t ,\ldots,
db^{*}_{1k}t ),\\
{\ve a}_2^{*}(t) & = ( db^{*}_{21}t,  db^{*}_{22}t+t_2 ,\ldots,
db^{*}_{2k}t ),\\
\vdots \\
{\ve a}_k^{*}(t) & = ( db^{*}_{k1}t,  db^{*}_{k2}t ,\ldots,
db^{*}_{kk}t+t_k ).\\
\end{array}$$

By Lemma \ref{Pol}, $\Lambda({\ve n}(t))$ is the polar lattice of
the lattice $\Lambda_{k+1}^{\bot}({\ve n}(t))$. Let $\{{\ve
a}_1(t),\ldots,{\ve a}_k(t)\}$ be a basis of $\Lambda({\ve n}(t))$
such that
$${\ve a}^{*}_i(t) {\ve a}_j(t)= \left \{
\begin{array}{ll}
1, & i=j,\\
0, & \mbox{otherwise}.\\
\end{array}
\right.$$ Consider the matrices
$A^{*}(t)=(a^{*}_{ij}(t))_{i,j=1}^{k}$ and
$B^{*}=(b^{*}_{ij})_{i,j=1}^{k}$. Let $A^{*}_{ij}(t)$ and
$B^{*}_{ij}$ be the minors obtained by omitting the $i$th row and
$j$th column in $A^{*}(t)$ and  $B^{*}$, respectively. Then, in
particular,

\begin{itemize}
\item[\hbox to 1.1cm{\rm (18)\hfill}]
$A^{*}_{ij}(t)=d^{k-1}t^{k-1}B^{*}_{ij}+O(t^{k-2}).$
\label{asympt_small_minors}\end{itemize} Moreover,
$${\ve
a}_i(t)=\lambda^{*}(A^{*}_{i1}(t),-A^{*}_{i2}(t)\ldots,(-1)^{k-1}A^{*}_{ik}(t)),$$
where $\lambda^{*}=\det \Lambda({\ve n}(t))=(\det
\Lambda_{k+1}^{\bot}({\ve n}(t)))^{-1}$. To check this, note that
$$\det \Lambda_{k+1}^{\bot}({\ve n}(t))={\ve
a}^{*}_i(t)(A^{*}_{i1}(t),-A^{*}_{i2}(t)\ldots,(-1)^{k-1}A^{*}_{ik}(t)).$$
Analogously, $${\ve
b}_i=\lambda(B^{*}_{i1},-B^{*}_{i2},\ldots,(-1)^{k-1}B^{*}_{ik}),$$
where $\lambda=(B^{*}_{k+1})^{-1}=(\det B^{*} )^{-1}$, since clearly
$$\det B^{*}={\ve
b}^{*}_i(B^{*}_{i1},-B^{*}_{i2}\ldots,(-1)^{k-1}B^{*}_{ik}).$$ By
(17)
$$\lambda^{*}= (d^{k}t^{k}\lambda^{-1}+O(t^{k-1}))^{-1}.$$

Thus by (18),
$$\begin{array}{l}
a_{ij}(t)=(-1)^{j-1}\frac{d^{k-1}t^{k-1}B^{*}_{ij}+O(t^{k-2})}{d^{k}t^{k}\lambda^{-1}+O(t^{k-1})}
=(-1)^{j-1}\frac{d^{k-1}t^{k-1}B^{*}_{ij}}{d^{k}t^{k}\lambda^{-1}(1+O\left(\frac{1}{t}\right))}+O\left(\frac{1}{t^2}\right)\\
= (-1)^{j-1}\frac{\lambda
B^{*}_{ij}}{dt}+O\left(\frac{1}{t^2}\right)
=\frac{b_{ij}}{dt}+O\left(\frac{1}{t^2}\right).\end{array}$$

\section{Proof of Lemma \ref{Super_Cool}}

The functions $f_t$, $f$ all are positive homogeneous of degree $1$.
Hence $f_t\rightarrow f$ uniformly on any bounded set. Thus

\begin{itemize}
\item[\hbox to 1.1cm{\rm (19)\hfill}]
${\ve l}_t\rightarrow {\ve l} \text{ implies } f_t({\ve
l}_t)\rightarrow f({\ve l}) \text{ as } t\rightarrow\infty .$
\label{L-1}
\end{itemize}

(i): Let $\epsilon>0$. Choose linearly independent vectors ${\ve
l}_1,\ldots,{\ve l}_k\in L$ such that

\begin{itemize}
\item[\hbox to 1.1cm{\rm (20)\hfill}]
$\max\{f({\ve l}_1),\ldots,f({\ve l}_i)\}\le \lambda_i(f,L)+\epsilon
\text{ for } i=1,\ldots,k.$ \label{L-2}
\end{itemize}
By Theorem 1 of \cite{GrLek}, pp.~178--179, there exist vectors
${\ve l}_{t1},\ldots,{\ve l}_{tk}\in L_t$ such that

\begin{itemize}
\item[\hbox to 1.1cm{\rm (21)\hfill}]
${\ve l}_{tj}\rightarrow {\ve l}_j \text{ as } t\rightarrow\infty
\text{ for } j=1,\ldots,k.$ \label{L-3}
\end{itemize}
Clearly,

\begin{itemize}
\item[\hbox to 1.1cm{\rm (22)\hfill}]
${\ve l}_{t1},\ldots,{\ve l}_{tk} \text{ are linearly independent
for all sufficiently large } t.$\label{L-4}
\end{itemize}
Thus
$$\begin{array}{l}
\lambda_i(f_t,L_t)\le\max\{f_t({\ve l}_{t1}),\ldots,f_t({\ve
l}_{ti})\} \le\max\{f({\ve l}_{1}),\ldots,f({\ve
l}_{i})\}+\epsilon\\
\le \lambda_i(f,L)+2\epsilon \text{ for } i=1,\ldots,k \text{ and
all sufficiently large } t
\end{array}$$
by (22), (21), (19) and (20), concluding the proof of (i).

(ii): Let $0<\epsilon<1$. Since $f_t\rightarrow f$ uniformly for
$||{\ve x}||=1$, $f({\ve x})>0$ for $||{\ve x}||=1$ and $f_t$ and
$f$ all are positive homogeneous of degree $1$, there is an
$\alpha>0$ such that

\begin{itemize}
\item[\hbox to 1.1cm{\rm (23)\hfill}]
$\alpha ||{\ve x}||\le(1-\epsilon)f({\ve x})\le f_t({\ve x}) \text{
for all } {\ve x} \text{ and all sufficiently large } t.$
\label{L-5}
\end{itemize}
For such $t$ the function $f_t({\ve x})$ is positive for ${\ve
x}\neq{\ve o}$. Thus the star body $\{{\ve x}: f_t({\ve x})\le 1\}$
is bounded. Hence we may choose

\begin{itemize}
\item[\hbox to 1.1cm{\rm (24)\hfill}]
${\ve l}_{t1},\ldots,{\ve l}_{tk}\in L_t, \text{ linearly
independent, such that}$\label{L-6}
\end{itemize}

$\max\{f_t({\ve l}_{t1}),\ldots,f_t({\ve
l}_{ti})\}=\lambda_i(f_t,L_t), i=1,\ldots,k \text{ for all
sufficiently large } t.$ By (23), (24) and (i),

\begin{itemize}
\item[\hbox to 1.1cm{\rm (25)\hfill}]
$\begin{array}{l} \displaystyle ||{\ve
l}_{tj}||\le\frac{1}{\alpha}\, f_t({\ve
l}_{tj})\le\frac{1}{\alpha}\,\lambda_i(f_t,L_t)\le
\frac{1}{\alpha}\,\lambda_d(f_t,L_t)\label{L-7}\\
\le\frac{1}{\alpha}\,\lambda_d(f,L)+\epsilon ,\ j=1,\ldots,k \text{
for all sufficiently large } t.\end{array}$
\end{itemize}
Moreover,
$$|\det({\ve l}_{t1},\ldots,{\ve l}_{tk})|\ge\det L_t\ge\det
L(1-\epsilon) \text{ for all sufficiently large } t$$ by (24) and
since $L_t\rightarrow L$ and thus $\det L_t\rightarrow \det L$.

The sequences $({\ve l}_{t1}), \ldots, ({\ve l}_{tk})$ all are
bounded by (25). The Bolzano -- Weierstrass theorem thus shows that
by considering suitable subsequences and re--indexing, if necessary,
we may assume that

$${\ve l}_{tj}\rightarrow {\ve l}_{j}\in L,\ |\det({\ve
l}_{1}, \ldots, {\ve l}_{k})|\ge \det L (1-\epsilon)>0\,,$$ see
\cite{GrLek}, pp.~178--179, Theorem~1. Hence ${\ve l}_{1}, \ldots,
{\ve l}_{k}$ are linearly independent and $f_t({\ve
l}_{tj})\rightarrow f({\ve l}_{j})$ by (19). Thus
$$\begin{array}{l}
\lambda_i(f_t, L_t)=\max\{f_t({\ve l}_{t1}),\ldots,f_t({\ve
l}_{ti})\}\rightarrow\max \{f({\ve l}_{1}),\ldots,f({\ve
l}_{i})\}\\
\ge \lambda_i(f,L),\ i=1,\ldots,k.\end{array}$$ Noting (i), this
concludes the proof of (ii).

\section{Proof of Theorem \ref{main_main}}

The inequality
$$C(K)=\sup_{
\begin{array}{c}
\scriptstyle{\ve n}\in\mathbb U^{k+1}
\end{array}}
\frac{\lambda_1(K,\Lambda({\ve n}))\cdots\lambda_k(K,\Lambda({\ve
n}))}{\det\Lambda({\ve n})}\le\frac{\alpha(K)}{\Delta(K)}$$ holds by
the definition of anomaly (see \cite{GrLek}, pp.~191, 192). To show
that equality holds, it is sufficient to prove that

\begin{itemize}
\item[\hbox to 1.1cm{\rm (26)\hfill}]
$\displaystyle\sup_{
\begin{array}{c}
\scriptstyle{\ve n}\in\mathbb U^{k+1}
\end{array}}
\frac{\lambda_1(K,\Lambda({\ve n}))\cdots\lambda_k(K,\Lambda({\ve
n}))}{\det\Lambda({\ve
n})}\ge\frac{\alpha(K)}{\Delta(K)}.$\label{let_us}
\end{itemize}
Let $\Lambda_0=\Lambda_0(K)$ be a lattice such that

\begin{itemize}
\item[\hbox to 1.1cm{\rm (27)\hfill}]
$\lambda_1(K,\Lambda_0)\cdots\lambda_k(K,\Lambda_0)=\frac{\alpha(K)}{\Delta(K)}\det\Lambda_0.$
\label{extreme_lattice}
\end{itemize}

The existence of such lattices for bounded star bodies in $\mathbb
R^{2}$ was proved in \cite{Mahler} and for all dimensions in
\cite{Hlawka}, see also \cite{Woods}. Let $\{{\ve r}_1,\ldots,{\ve
r}_k\}$ be a basis of $\Lambda_0$. Let $0<\delta<1$ and choose
linearly independent vectors ${\ve b}_1(\delta),\ldots,{\ve
b}_k(\delta)\in \mathbb Q^{k}$ such that

\begin{itemize}
\item[\hbox to 1.1cm{\rm (28)\hfill}]
$\begin{array}{l}||{\ve b}_j(\delta)-{\ve
r}_j||_\infty<\delta ,\ j=1,\ldots,k,\\
|\det({\ve b}^{T}_1(\delta),\ldots,{\ve
b}^{T}_k(\delta))-\det\Lambda_0|<\delta \det\Lambda_0.
\end{array}$ \label{diff_delta2}
\end{itemize}
Apply Theorem \ref{main_lemma} to the lattice $\Lambda$ with basis
$\{{\ve b}_1(\delta),\ldots,{\ve b}_k(\delta)\}$ and arbitrarily
chosen rational numbers $\alpha_1, \ldots,\alpha_k$ with
$0<\alpha_1\le\alpha_2\le\cdots\le\alpha_k\le 1$. This gives an
arithmetic progression $\cP$ and a sequence ${\ve
n}(t)=(n_1(t),\ldots,n_k(t),n_{k+1}(t))\in\mathbb Z^{k+1}$,
$t\in\cP$, such that $\Lambda({\ve n}(t))$ has a basis $\{{\ve
a}_1(t),\ldots,{\ve a}_k(t)\}$ where
$$dt a_{ij}(t)=b_{ij}(\delta)+O\left(\frac{1}{t}\right),\
i,j=1,\ldots,k.$$ Here $d=d(\delta)\in\mathbb N$ such that
$db_{ij}(\delta), d\alpha_jb_{ij}(\delta)\in\mathbb Z$ for all
$i,j=1,\ldots,k$. Choose any $t_0=t_0(\delta)\in\cP$ such that

\begin{itemize}
\item[\hbox to 1.1cm{\rm (29)\hfill}]
$||dt_0{\ve a}_j(\delta)-{\ve r}_j||_\infty<\delta,\ j=1,\ldots,k\,
\label{diff_delta}$
\end{itemize}
and $t_0>1/\delta$. Put $\Lambda_\delta=dt_0\Lambda({\ve n}(t_0))$.
For $\delta\rightarrow 0$ we obtain an infinite sequence of lattices
$\{\Lambda_\delta\}$ and by (29) $\Lambda_\delta\rightarrow
\Lambda_0$. In view of Corollary \ref{Cool},
$$\lambda_1(K,\Lambda_\delta)\cdots\lambda_k(K,\Lambda_\delta)\rightarrow\frac{\alpha(K)}{\Delta(K)}\det\Lambda_0\,
\text{ as } \delta\rightarrow 0.$$ We have
$$\lambda_1(K,\Lambda({\ve n}(t_0)))\cdots\lambda_k(K,\Lambda({\ve
n}(t_0)))=
\frac{\lambda_1(K,\Lambda_\delta)\cdots\lambda_k(K,\Lambda_\delta)}{(d(\delta)t_0(\delta))^{k}}
,$$ and by (3) and (28)
$$(d(\delta)t_0(\delta))^{k}= \frac{\det\Lambda}{\det\Lambda({\ve
n}(t_0))} + O(t_0^{k-1}) <
\frac{(1+\delta)\det\Lambda_0}{\det\Lambda({\ve n}(t_0))}(1 +
O(\delta )).$$ Thus, for every $\epsilon>0$ and for sufficiently
small $\delta>0$ there is an integer vector ${\ve n}={\ve
n}(t_0(\delta))$ such that

$$\lambda_1(K,\Lambda({\ve n}))\cdots\lambda_k(K,\Lambda({\ve
n}))>\frac{(1-\epsilon)\alpha(K)}{\Delta(K)}\det\Lambda({\ve n}).$$
This implies (26) and shows that (5) holds for the sequence $\{{\ve
n}(t_0(\delta))\}$. For this sequence equality (6) holds by (4) and
(7) holds by (3).

\section{Proof of Theorem \ref{Optimal}}

We shall show that for every $\epsilon>0$ there exists a vector
${\ve x}\in\mathbb R^{k}$ and a real number $Q>0$ such that

\begin{itemize}
\item[\hbox to 1.1cm{\rm (30)\hfill}]
$\{\lambda_1({\ve x},Q)\}^{k}>\frac{1-\epsilon}{\Delta(K)Q}.
\label{interval}$ \end{itemize} Let
$$C_1(K):=\limsup_{
\begin{array}{c}
\scriptstyle{\ve n}\in\mathbb U^{k+1}\\
\scriptstyle ||{\ve n}||_\infty\rightarrow\infty
\end{array}}
\frac{\{\lambda_1(K,\Lambda({\ve n}))\}^{k}}{\det\Lambda({\ve
n})}.$$ The proof of Theorem \ref{main_main} can be easily modified
to prove that

\begin{itemize}
\item[\hbox to 1.1cm{\rm (31)\hfill}]
$C_1(K)=\frac{1}{\Delta(K)}.\label{C_1}$
\end{itemize}
We just have to take for the lattice $\Lambda_0=\Lambda_0(K)$ any
critical lattice of $K$ and to replace (27) by the equality
$$\{\lambda_1(K,\Lambda_0)\}^{k}=\frac{\det\Lambda_0}{\Delta(K)}.$$

By (31) there is a sequence $\{{\ve n}(t)\}$, such that $||{\ve
n}(t)||_\infty\rightarrow\infty$ and for all sufficiently large $t$
holds
$$\{\lambda_1(K,\Lambda({\ve
n}(t)))\}^{k}>\frac{(1-\epsilon)\det\Lambda({\ve
n}(t))}{\Delta(K)\left(1-\frac{1}{n_{k+1}(t)}\right)}
=\frac{1-\epsilon}{\Delta(K)(n_{k+1}(t)-1)}.$$ Now put ${\ve
x}=(n_1(t)/n_{k+1}(t),\ldots,n_k(t)/n_{k+1}(t))$, $Q=n_{k+1}(t)-1$
and note that $\lambda_1(K,\Lambda({\ve n}(t)))=\lambda_1({\ve
x},Q)$.

\medskip
\noindent {\bf Remark.} The proof of Theorem \ref{Optimal} does not
yield only rational solutions ${\ve x}$ of the inequality (30) for
$\epsilon>0$. In fact, all vectors which are sufficiently close to a
vector ${\ve x}$ satisfying (30) satisfy (30) as well. Moreover,
since we apply Theorem \ref{main_lemma} with arbitrarily chosen
rational numbers $\alpha_i$, the equality (4) implies that solutions
of (30) approximate any rational point $(\alpha_1, \ldots,\alpha_k)$
with $0<\alpha_1\le\alpha_2\le\cdots\le\alpha_k\le 1$.

\section{Proof of Theorem \ref{th_on_decomp}}
\label{section_of_decompositions}

For any $\epsilon>0$ we have to find a sequence $\{{\ve n}(t)\}$ of
integer vectors such that $h({\ve n}(t))\rightarrow \infty$ and for
all sufficiently large $t$ the following inequality holds:

\begin{itemize}
\item[\hbox to 1.1cm{\rm (32)\hfill}]
$\displaystyle\inf_{
\begin{array}{c}
\scriptstyle{\ve p},\scriptstyle{\ve q}\in\mathbb Z^{k+1}\\
\scriptstyle\dim({\ve p},{\ve q})=2\\
\scriptstyle{\ve n}(t) = u{\ve p}+v{\ve q},u, v \in \mathbb Z
\end{array}
} \frac{h({\ve p})h({\ve q})}{h({\ve n}(t))^{1-\frac{1}{k}}}>
\frac{1-\epsilon}{(k+1)^{\frac{1}{k}}}. \label{S_e}$
\end{itemize}

Let ${\ve n} = (n_1, \ldots , n_{k+1})$, \mbox{$0< n_1\le \ldots\le
n_{k+1},$} be a primitive integer vector, that is \hbox{$\gcd(n_1,
\ldots , n_{k+1})=1$}, and let ${\ve m}=(m_1, \ldots, m_{k+1})$ be
an integer vector, such that ${\ve m}$ and ${\ve n}$ are linearly
independent. Consider the polygon $\Polygn=\Polygn({\ve m}, {\ve
n})$ defined by

\begin{itemize}
\item[\hbox to 1.1cm{\rm (33)\hfill}]
$\Polygn=\{ (x,y): |m_{i}y-n_{i}x|\le 1 \text{ for } i=1,\ldots ,
k+1\}. \label{Polygn}$
\end{itemize}
Let
\begin{itemize}
\item[\hbox to 1.1cm{\rm (34)\hfill}]
${\ve v}={\ve v}({\ve
m}):=\left(m_{1}-m_{k+1}\frac{n_1}{n_{k+1}},\ldots,m_{k}-m_{k+1}\frac{n_k}{n_{k+1}}\right)\in\Lambda({\ve
n}). \label{deltas}$
\end{itemize}
The following lemma is implicit in \cite{DecompII}.
\begin{lemma}
Let $0<n_1<\ldots < n_{k+1}$ and $\xi>0$. Then there is a centrally
symmetric convex set $ \cM_\xi = \cM_\xi({\ve n})\subset \mathbb
R^{k}$, such that ${\ve v}({\ve m})\in\cM_\xi$ for an integer vector
${\ve m}$ if and only if
$$\Delta(\Polygn({\ve m},
{\ve n}))\ge\frac{1}{n_{k+1}\xi}.$$ Moreover,

\begin{itemize}
\item[\hbox to 1.1cm{\rm (35)\hfill}]
$V_k(\cM_\xi)>(k+1)\xi^{k}.$ \label{M} \label{from_decompII}
\end{itemize}
\label{lemma_A}
\end{lemma}

Indeed, a set $\cM_\xi$ satisfying the equivalence stated in Lemma
\ref{lemma_A} is described by the formula (6) of \cite{DecompII} and
the inequality (35) is proved in Lemma~12 ibid. Let $f_{{\ve
n}}(\cdot)$ be the distance function of the set $\cM_1({\ve n})$. By
the definition of $\cM_\xi$, for ${\ve v}$ as in (34), we have that

\begin{itemize}
\item[\hbox to 1.1cm{\rm (36)\hfill}]
$f_{{\ve n}}({\ve v})=(n_{k+1}\Delta(\Polygn))^{-1}\,. \label{f_v}$
\end{itemize}
Consider a generalized honeycomb $E^{k}_1$ given by the inequalities
$$E^{k}_1=\{ x\in\mathbb R^{k}: |x_i|\le 1, |x_i-x_j|\le 1 \text{ for }\ i,j=1,\ldots,k,i\neq
j\}.$$ Observe that
$$E^{k}_1 =\bigcap_{p<q}\left\{ x\in \mathbb R^{k}: \left
(x_p, x_q \right ) \in E^2_1\right \}.$$ Let $g_k(\cdot)$ be the
distance function of $E^{k}_1$. Then clearly $$g_k({\ve
x})=\max_{1\le i<j\le k} g_2((x_i, x_j)).$$ By Lemma 1 of
\cite{DecompII},
$$V_k(E^{k}_1)=k+1, \Delta(E^{k}_1)=\frac{k+1}{2^{k}}$$
and $E^{k}_1$ has a unique critical lattice $\Lambda(E^{k}_1)$ with
basis
$$\begin{array}{llllllll}
{\ve b}_1 & =(  1, 1/2 ,\ldots, 1/2 ),\\
{\ve b}_2 & =( 1/2, 1 ,\ldots, 1/2 ),\\
\vdots\\
{\ve b}_k & =(  1/2, 1/2 ,\ldots, 1 ).
\end{array}$$

\begin{lemma}
For any $\epsilon>0$ there exists a $\delta=\delta(\epsilon)>0$ such
that for all integer vectors ${\ve n} = (n_1, \ldots ,n_k, n_{k+1})$
with $1-\delta<n_1/n_{k+1}<\ldots<n_k/n_{k+1}<1$, for all ${\ve
x}\in\mathbb R^{k}\setminus\{{\ve o}\}$
$$f_{{\ve n}}({\ve x})>\left(1-\frac{\epsilon}{2}\right)g_k({\ve
x}).$$ \label{f_g}
\end{lemma}

\begin{proof}
By formula (6) of \cite{DecompII}, the set $\cM_1({\ve n})$ is the
intersection of the sets $\cG_{pqr}$, where
$$\cG_{pqr}=\left\{x\in \mathbb R^{k}: (x_p, x_q)\in
\cB_1\left(\frac{n_p}{n_{k+1}}, \frac{n_q}{n_{k+1}}\right)\right\}$$
for $p<q<r=k+1$ and
$$\cG_{pqr}=\left \{x\in \mathbb R^{k}: \left
(x_p-\frac{n_p}{n_r}x_r, x_q-\frac{n_q}{n_r}x_r \right ) \in
\gamma\cB_1\left (\frac{n_p}{n_r}, \frac{n_q}{n_r} \right )\right
\}$$ for $p<q<r<k+1$, $\gamma=n_{k+1}/n_r$. The set
$B_1=B_1\left(\alpha, \beta\right)$, $0<\alpha<\beta<1$ is defined
by the formulae (8)--(13) of \cite{DecompI}. The boundary of $\cB_1$
consists of two horizontal segments
$$\pm
S_h=\left\{\pm( t,1)\in\mathbb R^{2}: -\frac{1-\alpha}{1+\beta}\le t
\le \frac{1+\alpha}{1+\beta}\right\}, \label{hor}$$ two vertical
segments
$$\pm S_v=\left\{\pm(1 , t)\in\mathbb R^{2}:
-\frac{1-\beta}{1+\alpha}\le t \le \frac{1-\beta}{1-\alpha}\right\},
\label{ver}$$ and four curvilinear arcs $\pm L_1$, $\pm L_2$ with
$$\begin{array}{l}
\displaystyle\pm L_1=\left\{\pm(x(t), tx(t))\in\mathbb R^{2}:
\frac{1-\beta}{1-\alpha}\le t\le
\frac{1+\beta}{1+\alpha}\right\},\\
x(t)=
\frac{-t^2(1+\alpha)^2+2t(1-\alpha+\beta+\alpha\beta)-(1-\beta)^2}
{4t(\beta-\alpha t)}\end{array}$$ and
$$\begin{array}{l} \pm
L_2=\left\{\pm(X(t), -tX(t))\in\mathbb R^{2}:
\frac{1-\beta}{1+\alpha}\le t\le
\frac{1+\beta}{1-\alpha}\right\},\\
X(t)=
\frac{-t^2(1-\alpha)^2+2t(1+\alpha+\beta-\alpha\beta)-(1-\beta)^2}
{4t(\beta+\alpha t)}.\end{array}$$ By Lemma 1 of \cite{DecompI},
$B_1$ is a centrally symmetric convex set.

Assume that there exists an $\epsilon>0$ such that for all
$\delta>0$, there exist an integer vector ${\ve n} = (n_1, \ldots ,
n_{k}, n_{k+1})$ with $1-\delta<n_1/n_{k+1}<\ldots<n_k/n_{k+1}<1$
and a point ${\ve x}\in\mathbb R^{k}\setminus\{{\ve o}\}$ with

\begin{itemize}
\item[\hbox to 1.1cm{\rm (37)\hfill}]
$f_{{\ve n}}({\ve x})\le\left(1-\frac{\epsilon}{2}\right)g_k({\ve
x}).\label{f_le_g}$
\end{itemize}
We shall show that this leads to a contradiction. By (37), there is
a point ${\ve a}=(a_1,\ldots,a_k)=\lambda {\ve x}$, $\lambda>0$,
such that $f_{{\ve n}}({\ve a})=1$ and

\begin{itemize}
\item[\hbox to 1.1cm{\rm (38)\hfill}]
$g_k({\ve a})=g_2((a_i,a_j))\ge
\left(1-\frac{\epsilon}{2}\right)^{-1} \label{outside}$
\end{itemize}
for some $i,j=1,\ldots,k$, $i<j$. Let $\alpha=n_i/n_{k+1}$,
$\beta=n_j/n_{k+1}$. Since ${\ve a}\in\cM_1({\ve n})$, we have
$(a_i, a_j)\in B_1(\alpha,\beta)$.

First, we consider the case $a_i a_j\ge 0$. By Lemma~2 of
\cite{DecompI}

\begin{itemize}
\item[\hbox to 1.1cm{\rm (39)\hfill}]
$B_1(\alpha,\beta)\subset C_1:=\{{\ve x}\in \mathbb R^{2}: ||{\ve
x}||_\infty\le 1\} \label{B_in_C}$
\end{itemize}
and thus
$$\{(x_i, x_j)\in B_1(\alpha,\beta): x_i x_j\ge
0\}\subset \{(x_i, x_j)\in E_1^2: x_i x_j\ge 0\},$$ which
contradicts (38).

Let us now consider the case $a_i a_j< 0$. Suppose $a_j=- t a_i$. We
may assume without loss of generality that

\begin{itemize}
\item[\hbox to 1.1cm{\rm (40)\hfill}]
$\left(1-\frac{\epsilon}{2}\right)^{-1}-1\le t \le
\left(\left(1-\frac{\epsilon}{2}\right)^{-1}-1\right)^{-1}
\label{t_in_interval}.$
\end{itemize}
Otherwise $(a_i,a_j)\not\in C_1$ and we get a contradiction with
(39). Since $(1-\beta)/(1+\alpha)$ tends to 0 and
$(1+\beta)/(1-\alpha)$ tends to infinity as $\delta$ tends to 0, we
have
$$\frac{1-\beta}{1+\alpha}< t<
\frac{1+\beta}{1-\alpha}$$ for $\delta$ small enough. Then $\mu(a_i,
a_j)\in\pm L_2$ for some $\mu\ge 1$. Further, for any $t$ from the
interval (40)
$$X(t)\rightarrow
\frac{1}{1+t},\ \text{ as } \delta\rightarrow 0.$$ Since
$g_2(1/(1+t), -t/(1+t))=1$, we obtain a contradiction with (38) for
all sufficiently small $\delta$.
\end{proof}

\begin{lemma}
For any $\epsilon>0$, there is an arithmetic progression $\cP$ and a
sequence of primitive integer vectors ${\ve
n}(t)=(n_1(t),\ldots,n_k(t),n_{k+1}(t))$, $t\in\cP$, such that
$h({\ve n}(t))\rightarrow \infty$ and for all sufficiently large
$t\in\cP$, for every non--zero vector ${\ve v}\in\Lambda({\ve
n}(t))$ the following holds
$$f_{{\ve n}(t)}({\ve
v})>(1-\epsilon)\left\{n_{k+1}(t)\Delta(E_1^{k})\right\}^{-\frac{1}{k}}.$$
\label{lemma_7_1}
\end{lemma}

\begin{proof}
Choose rational numbers
$1-\delta(\epsilon)<\alpha_1<\alpha_2<\cdots<\alpha_k< 1$ and apply
Theorem \ref{main_lemma} to the lattice $\Lambda=\Lambda(E^k_1)$,
the basis $\{{\ve b}_1,\ldots,{\ve b}_k\}$ of $\Lambda$ and the
numbers $\alpha_1,\alpha_2,\ldots,\alpha_k$. This yields an
arithmetic progression $\cP$ and a sequence of primitive integer
vectors ${\ve n}(t)$, $t\in\cP$ such that $h({\ve
n}(t))\rightarrow\infty$ and the corresponding lattices
$\Lambda({\ve n}(t))$ have bases ${\ve a}_1(t),\ldots,{\ve a}_k(t)$
where

\begin{itemize}
\item[\hbox to 1.1cm{\rm (41)\hfill}]
$\displaystyle
a_{ij}(t)=\frac{b_{ij}}{dt}+O\left(\frac{1}{t^2}\right),\
i,j=1,\ldots,k. \label{STAR}$
\end{itemize}
Here $d\in\mathbb N$ is such that $db_{ij},
d\alpha_jb_{ij}\in\mathbb Z$ for all $i,j=1,\ldots,k$. Moreover,
$$\alpha_i(t):=\frac{n_i(t)}{n_{k+1}(t)}
=\alpha_i+O\left(\frac{1}{t}\right).$$ Thus for sufficiently large
$t$,
$$1-\delta(\epsilon)<\frac{n_1(t)}{n_{k+1}(t)}<\ldots<\frac{n_k(t)}{n_{k+1}(t)}
<1.$$ We now show that for sufficiently large $t\in\cP$

\begin{itemize}
\item[\hbox to 1.1cm{\rm (42)\hfill}]
$\displaystyle \lambda_1(E_1^{k},\Lambda({\ve
n}(t)))>\left(1-\frac{\epsilon}{2}\right)\{n_{k+1}(t)\Delta(E_1^{k})\}^{-\frac{1}{k}}.$
\label{lambda_to_k}\end{itemize} The equality (41) implies that
$$dt\Lambda({\ve n}(t))\rightarrow\Lambda,
\text{ as } \ t\rightarrow\infty,\ t\in\cP.$$ Thus, Lemma
\ref{Super_Cool} implies that
$$\lambda_1(E_1^{k},dt\Lambda({\ve
n}(t)))\rightarrow 1, \text{ as } t\rightarrow \infty , t\in\cP .$$
Since
$$\lambda_1(E_1^{k},\Lambda({\ve n}(t)))=
\frac{\lambda_1(E_1^{k}, dt\Lambda({\ve n}(t)))}{dt}$$ and by (3),
$$dt=(n_{k+1}(t)\det\Lambda)^{\frac{1}{k}}\left(1 +
O\left(\frac{1}{t}\right)\right)^{\frac{1}{k}},$$ the inequality
(42) holds for all sufficiently large $t$. By Lemma~\ref{f_g} and
(42) for sufficiently large $t\in\cP$ for every non--zero vector
${\ve v}\in\Lambda({\ve n}(t))$,
$$\begin{array}{l}
f_{{\ve n}(t)}({\ve v})>\left(1-\frac{\epsilon}{2}\right)g_k({\ve
v})\ge\left(1-\frac{\epsilon}{2}\right)\lambda_1(E_1^{k},\Lambda({\ve
n}(t)))\\
>(1-\epsilon)\{n_{k+1}(t)\Delta(E_1^{k})\}^{-\frac{1}{k}}.\end{array}$$
The proof of the Lemma~\ref{lemma_7_1} is complete.
\end{proof}

After these preparations, the proof of Theorem~\ref{th_on_decomp} is
rather simple. We shall show that for every $\epsilon>0$ the
sequence $\{{\ve n}(t)\}_{t\in\cP}$ obtained in Lemma~5 satisfies
(32) for all sufficiently large $t$. Let $t\in\cP$ and let ${\ve p},
{\ve q}\in\mathbb Z^{k+1}$ be linearly independent vectors such that
${\ve n}(t)=u{\ve p}+v{\ve q}$ with $u,v\in\mathbb Z$, that is ${\ve
n}(t)\in\Lambda({\ve p}, {\ve q})$. Since the vector ${\ve n}(t)$ is
primitive, it can be extended to a basis of the lattice $S({\ve p},
{\ve q})\cap\mathbb Z^{k+1}$ by an integer vector ${\ve m}$.
Consider the polygon $\Polygn=\Polygn({\ve m}, {\ve n}(t))$ given by
(33). By Minkowski's lower bound for the product of successive
minima and since $V_2(\Polygn)\le 4\Delta(\Polygn)$, for all
linearly independent integer vectors $(x_1, y_1)$ and $(x_2, y_2)$,
$$\prod_{i=1}^{2}h(y_i{\ve m}-x_i{\ve
n}(t))\ge\lambda_1(\Polygn, \mathbb Z^{2})\lambda_2(\Polygn,\mathbb
Z^{2}) \ge2(V_2(\Polygn))^{-1}\ge
\frac{1}{2}(\Delta(\Polygn))^{-1}.$$ Since ${\ve p}, {\ve q}\in
\Lambda({\ve m}, {\ve n}(t))$, we have that

\begin{itemize}
\item[\hbox to 1.1cm{\rm (43)\hfill}]
$\displaystyle h({\ve p})h({\ve q})\ge
\frac{1}{2}(\Delta(\Polygn))^{-1}. \label{mult_lambda2}$
\end{itemize}
By (4), for all sufficiently large $t$ we have $h({\ve
n}(t))=n_{k+1}(t)$. Finally, by (43), (36) and
Lemma~~\ref{lemma_7_1}, for sufficiently large $t$, we get
$$\frac{h({\ve p})h({\ve q})}{h({\ve
n}(t))^{1-\frac{1}{k}}}\ge\frac{1}{2}
(n_{k+1}(t))^{\frac{1}{k}}f_{{\ve n}(t)}({\ve v}({\ve
m}))>\frac{1-\epsilon}{(k+1)^{\frac{1}{k}}}.$$

\vskip 1truecm

\noindent{\bf Acknowledgement.} The authors wish to thank Professor
A.~Schinzel for many valuable comments and suggestions. The first
author was supported by FWF Austrian Science Fund, projects M672,
M821--N12.

%
%

\end{document}